\documentstyle{amlts}
\begin{document}
\annalsline{156}{2002}
\received{October 12, 1999}
\startingpage{97}
\def\bye{\end{document}}
 \font\tenrm=cmr10
\def\binom#1#2{{#1\choose #2}}
\catcode`\@=11
\font\twelvemsb=msbm10 scaled 1100
\font\tenmsb=msbm10
\font\ninemsb=msbm10 scaled 800
\newfam\msbfam
\textfont\msbfam=\twelvemsb  \scriptfont\msbfam=\ninemsb
  \scriptscriptfont\msbfam=\ninemsb
\def\msb@{\hexnumber@\msbfam}
\def\Bbb{\relax\ifmmode\let\next\Bbb@\else
 \def\next{\errmessage{Use \string\Bbb\space only in math
mode}}\fi\next}
\def\Bbb@#1{{\Bbb@@{#1}}}
\def\Bbb@@#1{\fam\msbfam#1}
\catcode`\@=12

 \catcode`\@=11
\font\twelveeuf=eufm10 scaled 1100
\font\teneuf=eufm10
\font\nineeuf=eufm7 scaled 1100
\newfam\euffam
\textfont\euffam=\twelveeuf  \scriptfont\euffam=\teneuf
  \scriptscriptfont\euffam=\nineeuf
\def\euf@{\hexnumber@\euffam}
\def\frak{\relax\ifmmode\let\next\frak@\else
 \def\next{\errmessage{Use \string\frak\space only in math
mode}}\fi\next}
\def\frak@#1{{\frak@@{#1}}}
\def\frak@@#1{\fam\euffam#1}
\catcode`\@=12
\newcommand{\mgn}{\overline{M}_{g,n}}

\newcommand{\proj}{{\bf P}}
\newcommand{\barr}{\overline}
\newcommand{\rarr}{\rightarrow}
\newcommand{\oh}{{{\cal    O}}}
\newcommand{\com}{{\Bbb C}}
\newcommand{\Q}{{\Bbb Q}}
\newcommand{\Z}{{\Bbb Z}}
\newcommand{\G}{{\bf  G}}

\newcommand{\lan}{\langle}
\newcommand{\ran}{\rangle}
\newcommand{\eqq}{\stackrel{\sim}{=}}
\newcommand{\deli}{\bigtriangleup}
\newcommand{\sumo}{\oplus}
\newcommand{\grad}{\nabla}

\def\scup{{\hbox{\scriptsize$\cup$}}}
\def\scap{{\hbox{\scriptsize$\cap$}}}

 \newcommand{\LL}{{\Bbb L}}
\newcommand{\hodge}{{\Bbb E}}
\newcommand{\ch}{{\rm   ch}}
\newcommand{\rk}{{\rm   rk}}
\newcommand{\CM}{{\cal    M}}
\newcommand{\UM}{{\cal    C}}

\newcommand{\mata}{{\bf  A}}
\newcommand{\matb}{{\bf  B}}
\newcommand{\por}{\stackrel{\rarr}{{\cal    P}}}
\newcommand{\pord}{{\por}}
\newcommand{\ssym}{{\rm   sym}}
\newcommand{\matc}{{\bf  C}}
\newcommand{\matd}{{\bf  D}}

\font\eightmi=cmmi10


\title{Hodge integrals, partition matrices,\\ and the
%
$\lambda_g$ 
 conjecture}
\shorttitle{Hodge integrals} 

 \twoauthors{C. Faber}{R. Pandharipande}
 \institutions{Oklahoma State University, Stillwater, OK {\eightpoint and}\\
Institutionen f\"or Matematik, Kungl Tekniska H\"ogskolan, Stockholm, Sweden\\
{\eightpoint {\it E-mail addresses\/}: cffaber@math.okstate.edu, carel@math.kth.se}
\\ \vglue6pt
California Institute of Technology, Pasadena, CA, {\eightpoint and}\\
Princeton University, Princeton, NJ\\
{\eightpoint {\it E-mail addresses\/}: rahulp@cco.caltech.edu,
rahulp@math.princeton.edu}}

\centerline{\bf Abstract}
\vglue12pt

We prove a closed formula for integrals of the cotangent line 
classes against the top Chern class of the Hodge bundle on the 
moduli space of stable pointed curves. These integrals are computed via
       relations obtained from virtual localization in 
Gromov-Witten theory. An analysis of several natural matrices 
indexed by partitions is required.

\advance\sectioncount by -1

\section{Introduction}

0.1. {\it Overview}.
Let $M_{g,n}$ denote the moduli space of nonsingular genus $g$ curves with
$n$ distinct marked points (over $\com$).
Denote the 
moduli point corresponding the marked curve $(C,p_1,\ldots, p_n)$ by
$$[C,p_1, \ldots,p_n]\in M_{g,n}.$$ 
Let $\omega_C$ be the canonical bundle of algebraic
differentials on $C$.
The rank $g$ Hodge bundle,
$$\hodge \rarr M_{g,n},$$ has fiber $H^0(C, \omega_C)$ 
over $[C, p_1,\ldots,p_n]$.
The moduli space
$M_{g,n}$ is nonsingular of dimension 
$3g-3+n$ when considered as a stack (or orbifold). 

There is a natural compactification $M_{g,n} \subset \overline{M}_{g,n}$ by
stable curves (with nodal singularities).
The moduli space $\overline{M}_{g,n}$ is also a nonsingular stack.  
The Hodge bundle is well-defined
over $\overline{M}_{g,n}$: the fiber over a nodal curve $C$ is
defined to be the space of sections of the dualizing sheaf of $C$.
Let $\lambda_g $ be the top Chern class of $\hodge$ on
$\overline{M}_{g,n}$.
The main result of the paper is a formula for integrating tautological
classes on $\overline{M}_{g,n}$ against $\lambda_g$.

The study of integration against $\lambda_g$ has two main motivations.
First, such integrals arise naturally in the degree $0$ sector of the
Gromov-Witten theory of one-dimensional targets. The conjectural Virasoro 
constraints
of Gromov-Witten theory predict the $\lambda_g$ integrals have a
surprisingly simple form.
Second, the $\lambda_g$ integrals conjecturally govern the
entire tautological ring of the moduli space 
$$M_g^c\subset \overline{M}_g$$ of curves of compact type.
A stable curve is of compact type if the dual graph of $C$ is a tree.

\demo{{\rm 0.2.} Hodge integrals}
Let $A^*(\overline{M}_{g,n})$ denote the Chow ring of the moduli space
with $\Q$-coefficients.
We will consider two types of tautological classes
in $A^*(\overline{M}_{g,n})$:
\begin{itemize}
\item[$\bullet$] $\psi_i = c_1(\LL_i)$
for each marking $i$, where
$$\LL_i \rarr \overline{M}_{g,n}$$ denotes the cotangent line bundle
with fiber $T^*_{C,p_i}$ at the moduli point $[C,p_1,\ldots, p_n]\in \overline{M}_{g,n}$,
\item[$\bullet$] $\lambda_j = c_j(\hodge)$, for $j\leq g$.
\end{itemize}
Hodge integrals are defined to be the top intersection
products of the $\psi_i$ and $\lambda_j$ classes in 
$\overline{M}_{g,n}$. Hodge integrals 
play a basic role 
in Gromov-Witten theory and 
the study of the moduli space $\overline{M}_{g,n}$  
(see, for example, [Fa], [FaP1], [P]).
\enddemo

0.3. {\it Virasoro constraints and the $\lambda_g$ conjecture}.
The $\psi$ integrals in genus $0$ are determined by a
well-known formula:
\begin{equation}
\label{gen0}
\int_{\overline{M}_{0,n}} 
\psi_1^{\alpha_1} \cdots \psi_n^{\alpha_n} =
\binom{n-3}{\alpha_1,\ldots, \alpha_n}.
\end{equation}
The formula is a simple consequence of the
string equation [W].

The $\psi$ integrals are determined in all genera by
Witten's conjecture: the generating function of the
$\psi$ integrals satisfies the KdV hierarchy (or equivalently, 
Virasoro constraints).
Witten's conjecture has been proven by\break Kontsevich [K1].
A proof via Hodge integrals, Hurwitz numbers, and random trees
can be found in 
[OP].

The Virasoro constraints for the $\psi$ integrals over
$\overline{M}_{g,n}$ were
generalized to constrain tautological integrals
over the moduli space of stable maps to {\it     arbitrary}
nonsingular projective varieties through the work
of Eguchi, Hori, and Xiong [EHX], and Katz.
This generalization of Witten's original conjecture
remains open.

Tautological integrals over the moduli spaces of
{\it     constant} stable maps to nonsingular projective
varieties may be expressed as
Hodge integrals over $\overline{M}_{g,n}$. Hence,
the Virasoro constraints of [EHX] provide
(conjectural) constraints for Hodge integrals.
The $\lambda_g$ conjecture was found in [GeP] as
a consequence of these conjectural Virasoro constraints:
\begin{equation} \label{lamg}
\int_{\overline{M}_{g,n}} 
\psi_1^{\alpha_1} \cdots \psi_n^{\alpha_n} \lambda_g =
\binom{2g+n-3}{\alpha_1,\ldots, \alpha_n}\int_{{\overline{M}}_{g,1}}
\psi_1^{2g-2} \lambda_g,
\end{equation}
where $g\geq 1$, $\alpha_i\geq 0$.
In fact, conjecture (\ref{lamg}) was shown to be 
equivalent to the Virasoro constraints for constant maps
to an elliptic curve [GeP]. Equation~(\ref{lamg})
predicts
the combinatorics of the integrals
of the $\psi$ classes against
$\lambda_g$ is parallel to the genus $0$ formula (\ref{gen0}).
The integrals occurring in 
(\ref{lamg}) will be called $\lambda_g$ 
{\it     integrals}.

\demo{{\rm 0.4.} Moduli of curves of compact type}
The $\lambda_g$ integrals arise naturally in the
study of the moduli space of curves of compact type.
Let $$M_g^c \subset \overline{M}_g$$ denote the
(open) moduli space of curves of
compact type for $g\geq 2$.
The class $\lambda_g$ vanishes when
restricted to the complement $\overline{M}_g \setminus
M_g^c$ (see [FaP2]).
Integration against $\lambda_g$ therefore yields
a canonical linear evaluation function:
$$\epsilon: A^*(M^c_g) \rarr \Q,$$
$$\xi \in A^*(M^c_g), \ \ \epsilon(\xi) =
\int_{\overline{M}_g} \xi \cdot \lambda_g.$$
The $\lambda_g$ conjecture may be viewed as 
governing tautological evaluations in the 
Chow ring  $A^*(M^c_g)$.

The role of $\lambda_g$ in the study of $M_g^c$ exactly parallels
the role of $\lambda_g \lambda_{g-1}$ in the study of
$M_g$.
The class 
$\lambda_g \lambda_{g-1}$ vanishes
on the complement $\overline{M}_g \setminus M_g$. Hence, 
integration against $\lambda_g \lambda_{g-1}$ provides
a canonical evaluation function on $A^*(M_g)$ [Fa].

There is a conjectural formula for the $\lambda_g\lambda_{g-1}$
integrals which is also related to the Virasoro constraints [Fa], [GeP].
Data for $g\leq 15$ have led to a precise conjecture
for the ring of tautological classes $R^*(M_g) \subset A^*(M_g)$ [Fa].
In particular, $R^*(M_g)$ is conjectured to be Gorenstein
with the $\lambda_g\lambda_{g-1}$ integrals determining
the pairings into the socle.
It is natural to hope the tautological ring $R^*(M_g^c)\subset A^*(M_g^c)$
will also have a Gorenstein structure with socle pairings
determined by (\ref{lamg}).

A uniform perspective on the tautological rings $R^*(\overline{M}_g)$,
$R^*(M_g^c)$, and $R^*(M_g)$ 
may be found in [FaP2]. If the Gorenstein property holds
for $R^*(M_g^c)$, the $\lambda_g$ integrals determine the
entire ring structure [FaP2].
\enddemo

0.5. {\it Formulas for $\lambda_g$ integrals}.
The main result of the paper
is a proof of the $\lambda_g$ conjecture for all $g$.\pagebreak

\specialnumber{1}\proclaim{Theorem}
\label{mm}
The $\lambda_g$ integrals satisfy\/{\rm :}\/
$$\int_{\overline{M}_{g,n}} 
\psi_1^{\alpha_1} \cdots \psi_n^{\alpha_n} \lambda_g =
\binom{2g+n-3}{\alpha_1,\ldots, \alpha_n}\int_{{\overline{M}}_{g,1}}
\psi_1^{2g-2} \lambda_g\,.
$$
\endproclaim

The
integrals on the right side,
$$\int_{{\overline{M}}_{g,1}}
\psi_1^{2g-2} \lambda_g,$$
are determined by the following
formula previously proven in [FaP1]:

\begin{equation}
\label{sdsd}
F(t,k)=1+ \sum_{g\geq 1} \sum_{i=0}^{g} t^{2g} k^i
\int_{\overline{M}_{g,1}} \psi_1^{2g-2+i} \lambda_{g-i} = 
\Big( \frac{t/2}{\sin(t/2)} \Big)^{k+1}.\quad
\end{equation}
In particular, we find:
\begin{equation}
\label{qqq}
F(t,0)=
1+
\sum_{g\geq 1}t^{2g}\int_{{\overline{M}}_{g,1}}
\psi_1^{2g-2} \lambda_g = \Big( \frac{t/2}{\sin(t/2)} \Big).
\end{equation}
Equation (\ref{qqq}) is equivalent to the
Bernoulli number formula:
\begin{equation}
\label{fdfdfd}
\int_{{\overline{M}}_{g,1}} \psi_1^{2g-2}
\lambda_g  = \frac{2^{2g-1}-1}{2^{2g-1}} \frac{|B_{2g}|}{(2g)!}.
\end{equation}
Equation (\ref{fdfdfd}) and Theorem \ref{mm}
together determine all $\psi$ integrals against $\lambda_g$.

\demo{{\rm 0.6.} An interpretation in positive characteristic}
For an effective cycle $X$ on ${\overline{M}}_g$ with class
equal to a multiple of $\lambda_g$, the $\lambda_g$ conjecture may be viewed
as the analogue of Witten's conjecture for the family of curves
represented by~$X$. 

In characteristic $0$, it is not known whether
$\lambda_g$ is effective. In characteristic $p>0$ however, $\lambda_g$ {\it     is} effective.
Over an algebraically closed field of
characteristic $p$, define the $p$-rank $f(A)$ of an abelian variety by
$$p^{f(A)} = |A[p]|,$$
where $A[p]$ is the set of geometric $p$-torsion points.
Let $A_g$ be the moduli space
of principally polarized abelian varieties
of dimension $g$. Koblitz has shown the locus $V_0A_g$ of $p$-rank 0 abelian varieties
is complete and of codimension $g$ in $A_g$.
Van der Geer and Ekedahl [vdG] proved that the class of $V_0 A_g$ is 
proportional to $\lambda_g$
(by a factor equal to a polynomial in $p$).
Define the $p$-rank of a curve of compact type as the $p$-rank of its
Jacobian, and define the locus $V_0M_g^c$ of curves of $p$-rank $0$
via pullback along the Torelli morphism. This locus is complete in
${\overline{M}}_g$ and of \pagebreak codimension $g$ (see [FvdG]) --- it may however
be nonreduced. The  class of $V_0M_g^c$ is proportional
to $\lambda_g$ (by the same factor).
Hence $\lambda_g$ is effective in characteristic $p$.
The $\lambda_g$ conjecture may then be viewed as Witten's conjecture
for curves of $p$-rank~$0$. 

Perhaps this interpretation will eventually enhance our understanding
of the loci $V_0$. For example, $V_0A_g$ is expected to be
irreducible for $g\ge3$, but this is known only for $g=3$ (by a result of Oort).
The simple form of the Witten conjecture for $V_0M_g^c$ suggests
an analogy with genus $0$ curves that may lead to new insights.
\enddemo

0.7. {\it Localization}. Our proof of the
$\lambda_g$ conjecture uses the Hodge integral techniques 
introduced in [FaP1]. 
Let $\proj^1$ be equipped with an algebraic torus ${\bf  T}$ 
action. 
The virtual localization formula
established in [GrP] reduces all Gromov-Witten invariants 
(and their descendents) 
of $\proj^1$ to  explicit graph sums involving only 
Hodge integrals over $\mgn$.
Relations among the Hodge integrals may then be
found by computing invariants known to vanish.
The technique may be applied more generally by replacing
$\proj^1$ with any compact algebraic
homogeneous space. 

The philosophical basis of this method may be viewed
as follows.
If $M$ is an arbitrary
smooth variety with a torus action, the 
fixed
components of $M$ together with their equivariant normal
bundles satisfy global conditions obtained from the
geometry of $M$. Let $M$ be the (virtually) smooth moduli 
stack of
stable maps $\overline{M}_{g,n}(\proj^1)$ 
with the naturally induced  ${\bf  T}$-action. 
The ${\bf  T}$-fixed loci are then described as products of
moduli spaces of stable curves with virtual normal structures
involving the Hodge bundles [K2], [GrP]. 
In this manner, the geometry of $\overline{M}_{g,n}(\proj^1)$
imposes conditions on the ${\bf  T}$-fixed loci --- conditions
which may be formulated as relations among Hodge integrals
by [GrP].

Localization relations involving only the $\lambda_g$
integrals are found in Section 1 by studying
maps 
multiply
covering an exceptional $\proj^1$ of
an algebraic surface. 
These relations are
linear and involve a change of basis from the
standard form in formula (\ref{lamg}). 
However, it is not difficult to show the relations
are compatible with the $\lambda_g$ conjecture (see \S 2.4).
Both the linear equations from localization
and the change of basis are determined by natural
matrices indexed by partitions. In Section 3,
the ranks of these partition matrices are computed to
prove the system of linear equations found suffices
to determine all $\lambda_g$ integrals (up to the 
scalar $\int_{{\overline{M}}_{g,1}}
\psi_1^{2g-2} \lambda_g $ in each genus $g\geq 1$).

\demo{{\rm 0.8.} Acknowledgments}
We thank A.\ Buch, T.\ Graber, E.\ Looijenga, and R.\ Vakil
for several related conversations.
Discussions about partition matrices with D.\ Zagier
were very helpful to us.  
This project grew out of previous work with E.\ Getzler [GeP].
His ideas have played an important role in our research. 
The authors were partially supported by National Science
Foundation grants DMS-9801257 and DMS-9801574.
C.F.\ thanks the Max-Planck-Institut f\"ur Mathematik, Bonn,
for excellent working conditions and support, and the
California Institute of Technology for hospitality during
a visit in January/February 1999.\enddemo

\vglue-12pt

\section{\bf{Localization relations}}
\advance\eqcount by 5

1.1. {\it Torus actions}.
A system of linear
equations satisfied by the $\lambda_g$ integrals 
is obtained
here via localization relations. These relations are
found by computing vanishing integrals over moduli
spaces of stable maps in terms of Hodge integrals over
$\overline{M}_{g,n}$. 

The first step is to define the appropriate
torus actions.
Let $\proj^1=\proj(V)$ where $V=\com \oplus \com$.
Let $\com^*$ act diagonally on $V$:
\begin{equation}
\label{repp}
\xi\cdot (v_1,v_2) = ( v_1, 
\xi \cdot v_2).
\end{equation}
Let $p_1, p_2$ be the fixed points $[1,0], [0,1]$ of the corresponding
action on $\proj(V)$.
An equivariant lifting  of $\com^*$ to a line bundle $L$
over 
$\proj(V)$ is uniquely determined by the weights $[l_1,l_2]$
of the fiber
representations at the fixed points 
$$L_1= L|_{p_1}, \ \ \ L_2= L|_{p_2}.$$
The canonical lifting of $\com^*$ to the
tangent bundle $T_\proj$ has weights $[1,-1]$.
We will utilize the equivariant liftings of
$\com^*$ to $\oh_{\proj(V)}(1)$ and $\oh_{\proj(V)}(-1)$ with weights
$[1,0]$, $[0,1]$ respectively.
\vglue2pt

Let $\overline{M}_{g,n}(d)=\overline{M}_{g,n}(\proj(V), d)$ be the moduli
stack of stable genus $g$, degree $d$ maps to $\proj^1$ (see [K2], [FuP]).
There are canonical maps
$$\pi: U \rarr \overline{M}_{g,n}(d), 
\ \ \ \mu: U \rarr \proj(V)$$
where $U$ is the universal curve over the moduli stack.
The representation (\ref{repp}) canonically
induces $\com^*$-actions on $U$ and 
$\overline{M}_{g,n}(d)$ compatible
with the maps $\pi$ and $\mu$ (see [GrP]).

\demo{{\rm 1.2.} Equivariant cycle classes}
There are four  types of Chow classes in 
$A^*(\overline{M}_{g,n}(d))$ which will be considered here.
First,
there is a natural rank $d+g-1$ bundle on $\overline{M}_{g,n}(d)$: 
\begin{equation}
\label{wqwq}
{\Bbb R}=R^1\pi_* (\mu^* \oh_{\proj(V)}(-1)).
\end{equation}
The linearization $[0,1]$ on $\oh_{\proj(V)}(-1)$
defines an equivariant $\com^*$-action on ${\Bbb R}$.
Let $c_{\rm top}({\Bbb R})$ be the top Chern class
in $A^{g+d-1}(\overline{M}_{g,n}(d))$. Second,
the Hodge bundle
$$\hodge \rarr \overline{M}_{g,n}(d)$$
is defined by the vector space of differential
forms. There is a canonical lifting of the 
$\com^*$-action on $\overline{M}_{g,n}(d)$ to $\hodge$.
Let $\lambda_g \in A^{g}(\overline{M}_{g,n}(d))$ denote the top
Chern class of ${\Bbb E}$ as before.
Third,
for each marking $i$, let $\psi_i$
denote the first Chern class of 
the canonically linearized cotangent line corresponding to
$i$. Finally, let 
$${\rm   ev}_i: \overline{M}_{g,n}(d) \rarr \proj(V)$$
denote the $i^{\rm th}$ evaluation morphism, and let
$$\rho_i =c_1( {\rm   ev}_i^* \oh_{\proj(V)}(1)),$$
where we fix the $\com^*$-linearization $[1,0]$ on
$\oh_{\proj(V)}(1)$. \enddemo

1.3. {\it Vanishing integrals}.
A series of vanishing integrals $I(g,d,\alpha)$ 
over the moduli space of maps to $\proj^1$ 
is defined here. The parameters $g$ and $d$
correspond to the genus and degree of the map space.
Let $g\geq 1$
(the $g=0$ case is treated separately in \S 2.4).
Let
 $$\alpha=(\alpha_1, \ldots, \alpha_n)$$
be a (nonempty)
vector of nonnegative integers satisfying two conditions:
\begin{enumerate}
\item[(i)] $|\alpha|= \sum_{i=1}^n \alpha_i \leq d-2$,
\item[(ii)] $\alpha_i>0$ for $i>1$.
\end{enumerate}
By condition (i), $d\geq 2$.
Condition (ii) implies $\alpha_1$ is the only
integer permitted to vanish.
Let
\begin{equation}
I(g,d,\alpha)= 
\label{mvan1}
\int_{[\overline{M}_{g,n}(d)]^{\rm vir}}
\rho_1^{d-1-|\alpha|} \ \prod_{i=1}^n \rho_i \psi_i^{\alpha_i}
\ c_{\rm top}({\Bbb R}) \ \lambda_g.
\end{equation}
The virtual dimension of $\overline{M}_{g,n}(d)$
equals $2g+2d-2+n$. As the codimension of the
integrand equals $2g+2d-2+n$, the integrals are
well-defined. Since the class $\rho_1$ appears
in the integrand with exponent $d-|\alpha|\geq 2$ and
$\rho_1^2=0$, the integral vanishes.

These integrals occur in the following context.
Let $\proj^1 \subset S$ be an exceptional line
in a nonsingular algebraic surface. The virtual class
of the moduli space of stable maps to $S$ multiply covering
$\proj^1$ is obtained from the virtual class
of $\overline{M}_{g,n}(d)$ by intersecting with
$c_{\rm top}({\Bbb R})$. Hence, the series (\ref{mvan1})
may be viewed as vanishing Hodge integrals
over the moduli space of stable maps to~$S$.

\demo{{\rm 1.4.} Localization terms}
As all the integrand classes in the $I$ series
have been defined with $\com^*$-equivariant lifts,
the virtual localization formula of [GrP] yields
a computation of these integrals in terms
of Hodge integrals over moduli spaces of stable curves.

The integrals
(\ref{mvan1}) are
expressed as a  sum over connected decorated 
graphs~$\Gamma$ (see [K2], [GrP]) indexing the $\com^*$-fixed loci of 
$\overline{M}_{g,n}(d)$.
The vertices of these graphs lie over the
fixed points $p_1, p_2 \in \proj(V)$ and are
labelled with genera (which sum over the graph to $g-h^1(\Gamma)$).
The edges of the graphs lie over $\proj^1$ and
are labelled with degrees (which sum over the 
graph to $d$). Finally, the graphs carry
$n$ markings on the vertices. The edge valence of a vertex
is the number of incident edges (markings excluded).

In fact, only a very restricted subset of graphs 
will yield nonvanishing contributions to the $I$ series.
By our special choice of linearization
on the bundle ${\Bbb R}$, a vanishing result holds:
if a graph $\Gamma$
contains a vertex lying over $p_1$ of 
edge valence greater than 1, then the contribution of
$\Gamma$ to (\ref{mvan1}) vanishes.
A vertex over $p_1$ of edge valence at least 2
yields a trivial Chern root of ${\Bbb R}$ (with trivial weight $0$) in the
numerator of the localization formula to force the vanishing.
This basic vanishing was first used in $g=0$ by Manin in [Ma].
Additional applications have been pursued in [GrP], [FaP1].

By the above vanishing, only {\it     comb} graphs $\Gamma$ 
contribute to (\ref{mvan1}). Comb graphs contain
$k\leq d$ vertices lying over $p_1$ each
connected by a distinct edge to a unique vertex lying
over $p_2$. These graphs carry the usual vertex genus
and marking data.

Before deriving further restrictions on contributing
graphs, a classical result due to
Mumford is required [Mu].  \enddemo

\specialnumber{1}\proclaim{Lemma} Let $g\geq 1$. 
\label{mumm} 
$$ \sum_{i=0}^g \lambda_i  \cdot \sum_{i=0}^g (-1)^i \lambda_i =1$$
in $A^*(\overline{M}_{g,n})$.
In particular{\rm ,} $\lambda_g^2=0$.
\endproclaim

The factor $\lambda_g$ in the integrand of 
the $I$ series 
forces a further vanishing: if $\Gamma$ contains
a vertex over $p_1$ of positive genus, then the
contribution of $\Gamma$ to the integral
 (\ref{mvan1})  vanishes.
To see this, let $v$ be a positive genus $g(v)>0$
vertex lying over $p_1$. The integrand term $c_{\rm top}({\Bbb R})$
yields a factor $c_{g(v)}(\hodge^*)$ with trivial $\com^*$-weight on the
genus $g(v)$ moduli space corresponding to the vertex $v$.
The integrand class $\lambda_g$ factors as $\lambda_{g(v)}$
on each vertex moduli space. Hence, the equation
$$\lambda_{g(v)}^2=0$$ 
yields the required vanishing by Lemma \ref{mumm}.

The linearizations of the classes $\rho_i$ place restrictions on
the marking distribution.
As the class $\rho_i$ is
obtained from $\oh_{\proj(V)}(1)$ with linearization
$[1,0]$, all markings must lie on vertices over $p_1$
in order for the graph to contribute to (\ref{mvan1}).

Finally, we claim the markings of $\Gamma$ must lie on {distinct}
vertices over $p_1$ for nonvanishing contribution to
the $I$ series.
Let $v$ be a vertex over $p_1$ (with $g(v)=0$).
If $v$ carries at least two markings, the 
fixed locus corresponding to $\Gamma$ (see [K2], [GrP])
contains a product factor $\overline{M}_{0,m+1}$
where $m$ is the number of markings incident to $v$.
The classes $\psi_i^{\alpha_i}$
in the integrand of (\ref{mvan1}) carry trivial
$\com^*$-weight ---
they are pure Chow classes. Moreover, as each
$\alpha_i > 0$ for $i>1$, we see the sum of the
$\alpha_i$ as $i$ ranges over the set of markings
incident to $v$  is at least $m-1$. Since
this sum exceeds the dimension of $\overline{M}_{0,m+1}$,
the graph contribution to the $I$ series vanishes.

We have now proven the main result about the 
localization terms of the integrals (\ref{mvan1}).

\specialnumber{1}\proclaim{Proposition}
\label{qqqq}
The integrals in the $I$ series
are expressed via the virtual localization formula as
a sum over 
genus $g${\rm ,} degree $d${\rm ,} marked comb graphs $\Gamma$ satisfying\/{\rm :}
\begin{itemize}
\item[{\rm (i)}] all vertices over $p_1$ are of genus $0${\rm ,}
\item[{\rm (ii)}] 
each vertex over $p_1$ has at most one marking{\rm ,}
\item[{\rm (iii)}] 
the vertex over $p_2$ has no markings.
\end{itemize}
\endproclaim

1.5. {\it Hodge integrals}.
We introduce a new set of integrals over $\overline{M}_{g,n}$
which occur naturally in the
localization terms of  the $I$  series.
Let $g\geq 1$ 
(again the $g=0$ case is treated separately in \S  2.4).
Let $(d_1,\ldots ,d_k)$ be a nonempty sequence of positive integers.
Let
\begin{equation}
\label{fsfs}
\lan d_1,\ldots ,d_k \ran_g = \int_{\overline{M}_{g,k}} \frac{\lambda_g}
{\prod_{j=1}^k (1-d_j\psi_j)}.
\end{equation}
The value of the integral (\ref{fsfs})
clearly does not depend upon the ordering of
the sequence $(d_1, \ldots, d_k)$.

Let ${\cal    P}(d)$ denote the set of (unordered) partitions of $d>0$
into positive integers. Elements $P\in {\cal    P}(d)$
are unordered sets $P=\{d_1, \ldots, d_k\}$ of positive
integers with possible 
repetition.
The set ${\cal    P}(d)$ corresponds bijectively to the
set of distinct (up to reordering) {\it     degree $d$ integrals} by:
$$\{d_1, \ldots, d_k\} \mapsto  
\lan d_1, \ldots, d_k \ran_g 
$$
where  $\sum_{j=1}^k d_j =d$.

By the $\lambda_g$  conjecture,
we easily compute the prediction:
\begin{equation}
\label{preddy}
\lan d_1,\ldots, d_k \ran _g = \left( \sum_{j=1}^k d_j \right)^{2g-3+k}
 \int_{\overline{M}_{g,1}}
\psi_1^{2g-2} \lambda_g.
\end{equation}
\pagebreak

\noindent
Equation (\ref{preddy}) may be reduced further to the following
genus independent claim: for $g\geq 1$,
\begin{equation}
\label{xdxd}
\lan d_1,\ldots, d_k \ran _g = d ^{k-1} \lan
 d \ran_g
\end{equation}
where $\sum_{j=1}^k d_j=d$.
In Section 2.3, we will prove 
prediction (\ref{xdxd}) is {\it     equivalent} to the
$\lambda_g$ conjecture.

\demo{{\rm 1.6.} Formulas}
The precise contributions of allowable graphs
$\Gamma$ to the $I$ series
are now calculated.
Consider the integral $I(g,d, \alpha)$ where
$$\alpha=(\alpha_1, \ldots, \alpha_n).$$
Let $\Gamma$ be a
genus $g$, degree $d$, comb graph  with $n$ 
markings satisfying
conditions (i) and (ii) of Proposition \ref{qqqq}.
By condition (ii), $\Gamma$ must have  $k\geq n$ edges.
$\Gamma$ may be described uniquely
by the data 
\begin{equation}
\label{alll}
(d_1,\ldots, d_n)\  \scup \   \{ d_{n+1}, \ldots, d_k\},
\end{equation}
satisfying:
$$d_j >0,  \ \ \ \sum_{j=1}^k
d_j =d.$$
The elements of the ordered $n$-tuple $(d_1,\ldots, d_n)$
correspond
to the degree assignments of the edges incident
to the marked vertices. The elements of the unordered partition
$\{ d_{n+1}, \ldots, d_k\}$
correspond to the degrees of edges incident to
the unmarked vertices over $p_1$.
Let ${\rm   Aut}( 
\{ d_{n+1}, \ldots, d_k\})$ be the
group which permutes equal parts.
The group of graph 
automorphisms ${\rm   Aut}(\Gamma)$ (see [GrP]) 
equals ${\rm   Aut}(\{ d_{n+1}, \ldots, d_k\})$.

By a direct application of the virtual localization
formula of [GrP], we find the contribution of the graph
(\ref{alll})
to the (normalized) integral 
$$(-1)^{g+1}\cdot I(g,d, \alpha)$$
equals
$$\frac{1}{|{\rm   Aut}(\Gamma)|} \ 
\prod_{j=1}^n d_j^{-\alpha_j} \prod_{j=n+1}^{k} (-d_j)^{-1}
\prod_{j=1}^k \frac{d_j^{d_j}}{d_j!} \ \ \lan d_1,\ldots,d_k\ran_g.
$$
Hence, the vanishing of $I(g,d,\alpha)$ yields the Hodge
integral relation:
\begin{equation}
\label{fred}
\sum_{\Gamma} \frac{1}{|{\rm   Aut}(\Gamma)|} \ 
\prod_{j=1}^n d_j^{-\alpha_j} \prod_{j=n+1}^{k} (-d_j)^{-1}
\prod_{j=1}^k \frac{d_j^{d_j}}{d_j!} \ \ \lan d_1,\ldots,d_k \ran_g = 0,\hskip.4in
\end{equation}
where the sum is over all  graphs (\ref{alll}).

We point out two properties of the 
linear relations (\ref{fred}).
First, the relations do {\it     not}
depend upon the genus $g\geq 1$ ---
recall that the 
prediction (\ref{xdxd}) is also genus independent.
Second, the relations 
involve integrals
$\lan d_1, \ldots, d_k \ran_g$ with a fixed
sum $\sum_{j=1}^k d_j=d$.
By (\ref{fdfdfd}), the value $\lan d \ran_g$ is
never $0$.
Therefore, the integrals $\lan d_1, \ldots, d_k \ran_g$
are given at least one scalar dimension of freedom
in each degree $d$ by the equations (\ref{fred}). In Section
2.6, we will show that 
the solution space of the relations is exactly
one dimension in each degree.
\enddemo

1.7. {\it Generating functions}.
Let $g\geq 1$ as above.
Equation (\ref{fred})  may be rewritten
in a generating series form. While generating series will
not be used explicitly in our proof of Theorem \ref{mm}, the formalism
provides a concise description of the localization equations.

Let $t=\{t_1, t_2, t_3, \ldots\}$ be a set of variables indexed by
the natural numbers. Let $\Q[t]$ denote the
polynomial ring in these variables.
Define a $\Q$-linear function
$$\lan \ \ran : \Q[t] \rarr \Q$$
by the equations $\lan 1 \ran=1$ and
$$\lan t_{d_1} t_{d_2} \cdots t_{d_k} \ran = \lan
d_1, d_2, \ldots, d_k \ran _g.$$
We may extend $\lan \ \ran$ uniquely to define a $q$-linear function:
$$\lan \ \ran: \Q[t] [[q]] \rarr \Q[[q]].$$
For each nonnegative integer $i$, define:
$$Z_i(t,q) = \sum_{j>0} q^j t_j \frac{j^{j-i}}{j!} \in \Q[t][[q]].$$

The $I$ series equations (\ref{fred}) are equivalent to the
following constraints. 

\specialnumber{2}\proclaim{Proposition}
Let
$\alpha=(\alpha_1, \ldots, \alpha_n)$
be a  nonempty sequence of nonnegative
integers satisfying $\alpha_i>0$ for $i>1$. 
The series
$$\lan \exp(-Z_1) \cdot Z_{\alpha_1} \cdots Z_{\alpha_n} \ran
\in \Q[[q]]$$
is a {\rm     polynomial} of degree at most $1+\sum_{i=1}^n \alpha_i$
in $q$.
\endproclaim

\demo{Proof}
The coefficient terms of the expanded product 
$$\lan \exp(-Z_1) \cdot Z_{\alpha_1} \cdots Z_{\alpha_n} \ran$$
required to
vanish by the proposition coincide exactly with the relations 
(\ref{fred}). \phantom{lotsoffun}
\enddemo

1.8. {\it Example}.
Consider the polynomiality constraint obtained from the sequence
$\alpha=(0)$:
$${\rm   deg}_q\  \lan \exp(-Z_1) \cdot Z_0 \ran \ \ \leq 1.$$
After expanding the constraint, we find
$$\lan
\exp(-Z_1) \cdot Z_0 \ran
 = \lan t_1 \ran_g  q + (2 \lan t_2 \ran_g - \lan t_1^2 \ran_g ) 
q^2 +\cdots \ \ .$$ 
The equation 
\begin{equation}
\label{vvv}
2 \lan 2 \ran _g -  \lan 1,1 \ran_g =0
\end{equation}
is obtained from the $q^2$ term.
By the prediction (\ref{xdxd}), we see
equation (\ref{vvv}) is consistent with the
$\lambda_g$ conjecture.

\demo{{\rm 1.9.} The $\lambda_g$ conjecture}
The plan of the proof of the $\lambda_g$ conjecture
is as follows.
We first prove (\ref{xdxd}) is
equivalent to the $\lambda_g$ conjecture in Section 2.3.
The next step is to show the solution 
(\ref{xdxd}) satisfies all of our
linear relations (\ref{fred}). This result is
established in Section 2.4 via known $g=0$
formulas. In Section 2.6, 
the linear relations are proven to admit at most
a one-dimensional solution space in each degree.
Together, these three steps prove the 
$\lambda_g$ conjecture.

The above program relies upon the rank computations 
of certain natural matrices indexed by partitions.
The required results for these matrices are proven in
Section 3.\enddemo

\section{Proof of the $\lambda_g$ conjecture}

2.1. {\it String and dilaton}.
The $\lambda_g$ integrals 
satisfy the string and dilaton equations:
$$\int_{\overline{M}_{g,k+1}} \psi_1^{\alpha_1} \cdots \psi_k^{\alpha_k}
\psi_{k+1}^0
\lambda_g = 
\sum_{i=1}^k
\int_{\overline{M}_{g,k}} \psi_1^{\alpha_1} \cdots \psi_i^{\alpha_i-1}
\cdots \psi_k^{\alpha_k} \lambda_g,$$
$$\int_{\overline{M}_{g,k+1}} 
\psi_1^{\alpha_1} \cdots \psi_k^{\alpha_k}\psi_{k+1}^1
\lambda_g = 
(2g-2+k) \cdot
\int_{\overline{M}_{g,k}} \psi_1^{\alpha_1} 
\cdots \psi_k^{\alpha_k} \lambda_g.$$
The proofs of the string
and dilaton equations given in [W] are valid in the context 
of $\lambda_g$ integrals.

The
$\lambda_g$ conjecture is easily checked to be
compatible with 
the
string and dilaton
equations. For genus  $g=1$, all $\lambda_1$ integrals
must contain a $\psi_i^0$ 
factor in the integrand (for dimension reasons).
Hence, the $\lambda_1$ conjecture is
a consequence of the string equation.
Alternatively, the boundary equation
$12\lambda_1= \Delta_0$ in $A^1(\overline{M}_{1,1})$
immediately reduces the $\lambda_1$ conjecture to
the basic genus $0$ formula~(\ref{gen0}).

The $\lambda_g$
integrals for a fixed genus $g\geq 2$ may be expressed
in terms of
{\it     primitive integrals}
in which no factors $\psi_i^0$ or $\psi_i^1$
occur in the integrand. 
The distinct primitive integrals (up to ordering of the indices) are
in bijective correspondence with the
set ${\cal    P}(2g-3)$ of (unordered) partitions of $2g-3$.
The correspondence is given by
$$\{e_1, \ldots, e_l\}  \mapsto \int_{\overline{M}_{g,l}} 
\psi_1^{1+e_1} \ldots \psi_l^{1+e_l}
\lambda_g.$$
The value of the integral does not depend on the
ordering of the markings.
The $\lambda_g$ integrals may thus be viewed as having
${\cal    P}(2g-3)$ parameters in genus~$g$.

\advance\eqcount by 14
\demo{{\rm 2.2.} Matrix $\mata$}
Let $r\geq s >0$.
Let ${\pord}(r,s)$ be the set of
{\it     ordered} partitions of $r$ in exactly $s$
nonzero parts. An element $X \in \pord(r,s)$
is a vector $(x_1, \ldots, x_s)$.
Let $\mata$ be a matrix
with row and columns indexed by $\pord(r,s)$.
For $X,Y \in \pord(r,s)$,
define  the matrix element $\mata(X,Y)$ by:
$$\mata(X,Y)= \prod_{j=1}^s x_j^{-1+y_j}.$$

Let $\com^s$ be a vector space with coordinates
$z_1, \ldots, z_s$.
Let the partitions $X\in \pord(r,s)$
correspond to points in $\com^s$.
For $Y\in \pord(r,s)$, let $Y^-$ denote the
vector $(-1+y_1, \ldots, -1+y_s)$.
The set $\pord(r,s)$
corresponds bijectively
to the set of degree $r-s$ monomial functions in the
$z$ variables by:
\begin{equation}
\label{bbnn}
Y  \ \ \leftrightarrow \ \ m_{Y^-}(z)=z_1^{-1+y_1} \cdots z_s^{-1+y_s}.
\end{equation}
$\mata$ is simply the matrix obtained by evaluating
degree $r-s$ monomials on partition
points in $\com^s$:
$$\mata(X,Y)= m_{Y^-}(x_1, \ldots, x_s).$$
The following lemma needed here
will be proven in Section 3.1.

\specialnumber{2}\proclaim{Lemma} 
\label{lmaa}
For all pairs $(r,s)${\rm ,} the matrix
$\mata$ is invertible.
\endproclaim

The symmetric group ${\Bbb S}_s$ acts naturally
by permutation on the set $\pord(r,s)$.
Let $V_{r,s}$ denote the canonically induced
${\Bbb S}_s$ permutation representation.
The matrix $\mata$ determines a natural
${\Bbb S}_s$-invariant
bilinear form:
$$\phi: V_{r,s} \times V_{r,s} \rarr \com$$
by $\phi([X],[Y]) = \mata(X,Y)$.
The form $\phi$ is nondegenerate by Lemma \ref{lmaa}.
Let $V_{r,s}^{{\Bbb S}}\subset V_{r,s}$ denote the
${\Bbb S}_s$ invariant subspace.
By an application of Schur's Lemma, the restricted form
$$\phi^{{\Bbb S}}: V^{{\Bbb S}}_{r,s} \times 
V^{{\Bbb S}}_{r,s} \rarr \com$$
is also nondegenerate.

Let ${\cal    P}(r,s)$ denote the
set of (unordered) partitions of $r$ in exactly
$s$ parts.
An element $P\in {\cal    P}(r,s)$
is a set $\{ p_1, \ldots, p_s\}$ of positive integers (with possible
repetition).
The set ${\cal    P}(r,s)$
may be placed in bijective correspondence with
a basis of $V_{r,s}^{{\Bbb S}}$ by
\begin{equation}
\label{vvvv}
\{ p_1, \ldots, p_s\} \ \  \leftrightarrow \ \ \sum_{\sigma \in {\Bbb S}_s}
[(p_{\sigma(1)}, \ldots, p_{\sigma(s)})].
\end{equation}

The correspondence (\ref{bbnn}) yields
an equivariant isomorphism between $V_{r,s}$ and
the vector space of polynomial functions of homogeneous
degree $r-s$ in the $z$ variables.
Via this isomorphism,
the basis element (\ref{vvvv}) corresponds to the
symmetric function:
$$\ssym(m_{P^-})= \sum_{\sigma\in {\Bbb S}_s} z_1^{-1+p_{\sigma(1)}}
\cdots z_s^{-1+p_{\sigma(s)}}.$$

In the basis (\ref{vvvv}), the form $\phi^{{\Bbb S}}$ corresponds
to the matrix
$\mata^{{\Bbb S}}$ with rows and columns
indexed by ${\cal    P}(r,s)$
and matrix element
$$\mata^{{\Bbb S}}(P,Q)= {s!} \cdot
 \ssym(m_{Q^-}) \Big(p_1,\ldots,p_s \Big).$$
As a corollary of Lemma \ref{lmaa}, we have proven:

\specialnumber{3}\proclaim{Lemma} 
\label{lmaas}
For all pairs $(r,s)${\rm ,} the matrix
$\mata^{{\Bbb S}}$ is invertible.
\endproclaim

2.3. {\it Change of basis}.
The partition matrix results
of 
Section 2.2 are required for the
following proposition. This is the first step
in the proof of the $\lambda_g$ conjecture.

\specialnumber{3}\proclaim{Proposition}
\label{cbas}
Let $g\geq 2$.
The values of the primitive $\lambda_g$ integrals are
uniquely determined by the degree $2g-3$ integrals\/{\rm :}
$$\{ \lan d_1, \ldots, d_k \ran_g \}$$
where $\sum_{j=1}^k d_j=2g-3$.
\endproclaim

\demo{Proof}
Let $D=\{d_1, \ldots, d_k\} \in {\cal    P}(2g-3,k)$.
We may certainly express the integral
$$\lan D \ran _g =\lan d_1, \ldots, d_k \ran_g$$
in terms of the primitive $\lambda_g$ integrals by:
\begin{equation}
\label{bnm}
\lan D \ran_g = \sum_{l=1}^{k}\ 
\sum_{E=\{e_1,\ldots,e_l\}
\in {\cal    P}(2g-3,l)} M(D,E) \cdot \int_{\overline{M}_{g,l}} 
\psi_1^{1+e_1} \ldots \psi_l^{1+e_l}
\lambda_g. \hskip.25in
\end{equation}
Note  no primitive $\lambda_g$ integrals corresponding
to partitions of length greater than $k$ occur in the sum.
The string and dilaton equations are required  
to compute the values $M(D,E)$ 
where the length of $E$ is strictly less than $k$.

Let
$M$ be the matrix with
rows and columns indexed by 
${\cal    P}(2g-3)$ and matrix elements $M(D,E)$.
In order to establish the proposition, it suffices to \pagebreak
prove $M$ is invertible.

We order the rows and columns of $M$ by increasing length of
partition (the order within a fixed length can be
chosen arbitrarily).
$M$ is then block lower-triangular with diagonal blocks $M_k$ 
determined by 
partitions of a fixed length $k$. Hence,
$$\det(M)= \prod_{k=1}^{2g-3} \det(M_k).$$
We will
prove $\det(M_k) \neq 0$ for each $k$.

Let $k$ be a fixed length.
The diagonal block $M_k$ 
has rows and columns indexed by ${\cal    P}(2g-3,k)$.
Let $D, E \in {\cal    P}(2g-3,k)$.
The matrix element $M_k(D,E)$ is
given by:
$$M_k(D,E) = \frac{ \sum_{\sigma \in 
{\Bbb S}_k}  \prod_{j=1}^k d_j^{1+e_{\sigma(j)}}}
{|{\rm   Aut}(E)|} =
\frac{\prod_{j=1}^k d_j^2}{|{\rm   Aut}(E)|} \cdot 
\ssym(m_{E^-}) \Big(d_1, \ldots, d_k \Big)
.$$
Here ${\rm   Aut}(E)$ is the group permuting equal parts of 
the partition $E$.
This element is computed by a simple expansion of the
denominator in the definition~(\ref{fsfs}) of 
the integral $\lan D \ran_g$. No applications of the
string or dilaton equations are necessary.

Let
$\mata^{{\Bbb S}}$ be the matrix
defined in Section 2.2 for $(r,s)=(2g-3,k)$.
For $X,Y \in {\cal    P}(2g-3,k)$,
$$\mata^{{\Bbb S}} (P,Q) = {k!} \cdot
\ssym(m_{Q^-}) \Big(p_1, \ldots, p_k \Big)
.$$
As $M_k$ differs from $\mata^{{\Bbb S}}$ 
only by scalar row and column operations, $M_k$ 
is invertible if and only if
$\mata^{{\Bbb S}}$ is invertible. However, 
by Lemma \ref{lmaas}, 
$\mata^{{\Bbb S}}$ is invertible. \phantom{yetmorefun}
\enddemo

By Proposition \ref{cbas}, the $\lambda_g$ conjecture
is equivalent to the prediction:
 \begin{equation}
\label{ydyd}
\lan d_1,\ldots, d_k \ran _g = d ^{k-1} \lan
 d \ran_g
\end{equation}
where $\sum_{j=1}^k d_j=d$.
We will prove the $\lambda_g$ conjecture in 
form (\ref{ydyd}).

\demo{{\rm 2.4.} Compatibility}
We now prove equation (\ref{ydyd}) yields a solution of the
linear system of equations obtained from
localization (\ref{fred}).
Our method is to use localization equations
in genus $0$ together with the basic formula (\ref{gen0}).

Define genus $0$ integrals $\lan d_1, \ldots, d_k \ran_0$
by:
\begin{equation}
\label{fsfsz}
\lan d_1,\ldots ,d_k \ran_0 = \int_{\overline{M}_{0,k+2}} \frac{\lambda_0}
{\prod_{j=1}^k (1-d_j\psi_j)} =\int_{\overline{M}_{0,k+2}} \frac{1}
{\prod_{j=1}^k (1-d_j\psi_j)}. \hskip.27in
\end{equation}
As $k+2\geq 3$, these integrals are well-defined (the
two extra markings of $\overline{M}_{0,k+2}$ serve
to avoid the degenerate spaces $\overline{M}_{0,1}$ and
$\overline{M}_{0,2}$).
An easy evaluation using (\ref{gen0}) shows:
\begin{equation}
\label{wqw}
\lan d_1, \ldots, d_k \ran_0 = d^{k-1}, \ \ \ \sum_{j=1}^k d_j =d.
\end{equation}
In particular,
$$\lan d_1, \ldots, d_k \ran_0 = d^{k-1} \lan d \ran_0.
$$

Relations among the integrals $\lan d_1,\ldots, d_k \ran_0$
may be found in a manner similar to the
higher genus development in Section 1.
We follow the notation of the $\com^*$-action on
$\proj^1$ introduced in Sections 1.1--1.2.
The
$\com^*$-equivariant classes 
$$c_{\rm top}({\Bbb R}), \psi_i, \rho_i$$
are defined on the moduli space $\overline{M}_{0,n}(d)$.
Define a new class
$$\gamma_i = c_1( {\rm   ev}_i^*( \oh_{\proj(V)}(-1) ))$$
with $\com^*$-linearization determined by the
action with weights $[0,1]$ on the line bundle
$\oh_{\proj(V)}(-1)$.

Again, we find a series $I(0,d,\alpha)$ of vanishing integrals.
We require $\alpha$ to satisfy conditions (i) and (ii)
of Section 1.3.
\begin{equation}
I(0,d,\alpha)= 
\label{mvan10}
\int_{[\overline{M}_{0,n+2}(d)] }
\rho_1^{d-1-|\alpha|} \ \prod_{i=1}^n \rho_i \psi_i^{\alpha_i}
\ c_{\rm top}({\Bbb R}) \ \gamma_{n+1} \gamma_{n+2}.
\end{equation}
These integrals 
are well-defined and vanish as before.

The localization formula yields a computation of the
vanishing integrals~(\ref{mvan10}). The argument
exactly follows the higher genus development in
Sections 1. In addition to the graph restrictions
found in Section 1.4,
the two extra points (corresponding to the $\gamma$ factors
in the integrands) must lie on the
unique vertex over the fixed point $p_2 \in \proj(V)$.
These extra points ensure that the unique vertex over $p_2$
will not degenerate in the localization formulas. The
resulting graph contributions then agree exactly with the
expressions found in Section 1.6.

$I(0,d,\alpha)$ yields the 
relation:
\begin{equation}
\label{fredd}
\sum_{\Gamma} \frac{1}{|{\rm   Aut}(\Gamma)|} \ 
\prod_{j=1}^n d_j^{-\alpha_j} \prod_{j=n+1}^{k} (-d_j)^{-1}
\prod_{j=1}^k \frac{d_j^{d_j}}{d_j!} \ \ \lan d_1,\ldots,d_k \ran_0 = 0, \hskip.35in
\end{equation}
where the sum is over all graphs: 
$$\Gamma = (d_1,\ldots, d_{n}) \  \scup \    \{ d_{n+1}, \ldots, d_k\}, \ \ \
 d_j>0, \ \ \ \sum_{j=1}^k d_j = d.
$$

Equation (\ref{fredd}) equals the specialization of
equation (\ref{fred}) to genus $0$. 
Hence, we have proven
the predicted form proportional to (\ref{wqw}) solves the linear relations
obtained from \pagebreak localization. 
\enddemo

2.5. {\it Matrix $\matb$}.
Let $r>s >0$. As in Section 2.2, let
$\com^s$ be a vector space with coordinates $z_1, \ldots, z_s$.
Let the set $\pord(r,s)$ correspond
to points in $\com^s$ by the new association:
$$X \in  \pord(r,s)\ \  \leftrightarrow \ \ 
\left(\frac{1}{x_1}, \ldots,\frac{1}{x_s}\right)
\in \com^s.$$
Let ${\cal    M}(r,s)$ be the set of 
monomials $m(z)$
in the
coordinate variables satisfying the following two conditions:
\begin{enumerate}
\item[(i)] ${\rm   deg}(m) \leq r-2$,
\item[(ii)] $m(z)$ omits {\it     at most} one coordinate factor $z_i$.
\end{enumerate}
Note the condition ${\rm   deg}(m) \geq s-1$ is a consequence of
condition (ii).
The set ${\cal    M}(r,s)$ is never empty. 

Let $\matb$ be a matrix with rows indexed by
${\cal    M}(r,s)$ and columns indexed by
$\pord (r,s)$.
Let the matrix element $\matb(m,X)$ be defined by
evaluation:
$$ \matb(m,X)= m\left(\frac{1}{x_1},\ldots, 
\frac{1}{x_s}\right).$$ The following lemma will be
proven in Section 3.2.

\specialnumber{4}\proclaim{Lemma} 
\label{bbbb}
For all pairs $(r,s)${\rm ,} the matrix $\matb$
has rank equal to $|\pord(r,s)|$.
\endproclaim

There is a natural ${\Bbb S}_s$-action on the
set ${\cal    M}(r,s)$ defined by:
\begin{equation}
\label{acty}
\sigma(z_1^{\alpha_1} \cdots z_s^{\alpha_s})= 
z_1^{\alpha_{\sigma(1)}} \cdots z_s^{\alpha_{\sigma(s)}}.
\end{equation}
Let $W_{r,s}$ denote the ${\Bbb S}_s$
permutation representation induced by the action (\ref{acty}).
As before, let $V_{r,s}$ denote the
${\Bbb S}_s$ permutation representation
induced by the natural group action on $\pord(r,s)$.

The matrix $\matb$ determines a natural
${\Bbb S}_s$-invariant
bilinear form:
$$\phi: W_{r,s} \times V_{r,s} \rarr \com$$
by $\phi([m],[X]) = \matb(m,X)$.
The form $\phi$ induces 
a canonical homomorphism of ${\Bbb S}_s$ representations:
$$W_{r,s} \rarr V_{r,s}^* \rarr 0,$$
surjective by Lemma \ref{bbbb}.
By Schur's lemma, the restricted morphism is also surjective:
$$W^{{\Bbb S}}_{r,s} \rarr V^{{\Bbb S}*}_{r,s} \rarr 0.$$
Hence the restricted form:
$$\phi^{{\Bbb S}}: W^{{\Bbb S}}_{r,s} 
\times V^{{\Bbb S}}_{r,s} \rarr \com$$
has rank equal to $|{\cal    P}(r,s)|$.

Let ${\cal    M}_{\ssym}(r,s)$ denote the
set of distinct symmetric functions obtained
by symmetrizing monomials in ${\cal    M}(r,s)$:
$$m \in {\cal    M}(r,s) \ \ \rarr \ \ {\rm   sym}(m)=
\sum_{\sigma\in {\Bbb S}_s} 
\sigma(m).$$
The set ${\cal    M}_{\ssym}(r,s)$ corresponds to a basis of
$W^{{\Bbb S}}_{r,s}$.
Let the set ${\cal    P}(r,s)$ correspond to a basis
of $V^{{\Bbb S}}_{r,s}$ as before (\ref{vvvv}).

Let $\matb^{{\Bbb S}}$  
be a matrix  with rows indexed by ${\cal    M}_{\ssym}(r,s)$, 
columns indexed by ${\cal    P}(r,s)$,
and matrix element:
$$\matb^{{\Bbb S}}({\rm   sym}(m),P)= {s!}
 \cdot \ssym(m) \Big(\frac{1}{p_1}, \ldots, \frac{1}{p_s} \Big).$$
The restricted form $\phi^{{\Bbb S}}$
expressed in the bases ${\cal    M}_{\ssym}(r,s)$ and
${\cal    P}(r,s)$  corresponds to the
matrix $\matb^{{\Bbb S}}$.
As a corollary of Lemma \ref{bbbb}, we have proven:

\specialnumber{5}\proclaim{Lemma} 
\label{lmbbs}
For all pairs $(r,s)${\rm ,} the matrix
$\matb^{{\Bbb S}}$ has rank equal to $|{\cal    P}(r,s)|$.
\endproclaim

2.6. {\it Linear relations}.
The rank computation of $\matb^{{\Bbb S}}$ directly
yields the final step in the proof of the $\lambda_g$
conjecture.

\specialnumber{4}\proclaim{Proposition}
\label{ffff}
Let $d\geq 1$.
The linear relations {\rm (\ref{fred})}
admit at most a one\/{\rm -}\/dimensional solution space
for the integrals
\begin{equation}
\label{zfg}
\lan d_1, \ldots, d_k\ran_g, \  \ \ \ \sum_{j=1}^k d_j =d.
\end{equation}
\endproclaim

\demo{Proof}
As
no linear relations in (\ref{fred})
constrain the unique degree 1 integral $\lan 1 \ran_g$,
we may assume $d\geq 2$.

Recall the distinct integrals (\ref{zfg}) 
correspond to the set ${\cal    P}(d)$. 
There is a unique integral of partition length $d$:
$$\lan 1,\ldots,1 \ran_g.
$$
We will prove that the localization relations determine
all degree $d$ integrals in terms of $\lan 1,\ldots,1 \ran_g$.

We proceed by descending induction on the partition length.
If $D\in {\cal    P}(d)$ is of length $l(D)=d$, then
$\lan D \ran_g$ equals $\lan 1, \ldots,1 \ran_g$
 --- the base
case of the induction.

Let $d>n>0$.
Assume now all integrals corresponding to partitions 
$D\in {\cal    P}(d)$
of
length {\it     greater} than $n$ are determined in terms
of $\lan 1, \ldots,1 \ran_g$.
 Consider the integrals corresponding
to the partitions ${\cal    P}(d,n)$.
For each nonempty sequence
$$\alpha=(\alpha_1, \ldots, \alpha_n)$$
satisfying
\begin{itemize}
\item[(i)] $|\alpha|= \sum_{i=1}^n \alpha_i \leq d-2$,
\item[(ii)] $\alpha_i>0$ for $i>1$,
\end{itemize}
we obtain the relation:
\begin{equation}
\label{fredq}
\sum_{\Gamma} \frac{1}{|{\rm   Aut}(\Gamma)|} \ 
\prod_{j=1}^n d_j^{-\alpha_j} \prod_{j=n+1}^{k} (-d_j)^{-1}
\prod_{j=1}^k \frac{d_j^{d_j}}{d_j!} \ \ \lan d_1,\ldots,d_k \ran_g = 0.\hskip.4in
\end{equation}
Recall the sum is over all graphs:
$$\Gamma = (d_1,\ldots, d_{n}) \  \scup \    \{ d_{n+1}, \ldots, d_k\}, \ \ \
 d_j>0, \ \ \ \sum_{j=1}^k d_j = d.
$$

We note only integrals corresponding to partitions of
length at least $n$ occur in (\ref{fredq}).
By the inductive assumption,
only the terms in (\ref{fredq}) containing integrals
of length {\it     exactly} $n$ concern us:
\begin{equation}
\label{fredqq}
\sum_{\Gamma}  
\prod_{j=1}^n d_j^{-\alpha_j} 
\prod_{j=1}^n \frac{d_j^{d_j}}{d_j!} \ \ \lan d_1,\ldots,d_n \ran_g = 
f_{\alpha}( \lan 1,\ldots,1 \ran_g).
\end{equation}
The sum is over all ordered sequences:
$$(d_1, \ldots, d_n), \ \ \ d_j>0, \ \ \ \sum_{j=1}^n d_j=d.$$
The factor $|{\rm   Aut}(\Gamma)|$ is trivial for the
terms containing integrals of length exactly $n$.

Let $L_\alpha$ denote the linear equation (\ref{fredqq}).
To each $\alpha$, we may associate
an element of
${\cal    M}(d,n)$ by
$$\alpha \rarr m_\alpha = z_1^{\alpha_1} \ldots z_n^{\alpha_n}.$$
Let $D\in {\cal    P}(d,n)$. 
The coefficient of $\lan D \ran_g$ in $L_\alpha$ is
$$\frac{  \prod_{j=1}^n \frac{d_j^{d_j}}{d_j!}  }{|{\rm   Aut}(D)|} \ 
{\rm   sym}(m_{\alpha})
\Big( \frac{1}{d_1},\ldots,\frac{1}{d_n} \Big).$$
As before, ${\rm   Aut}(D)$ is the group permuting
equal parts of $D$. The equation $L_\alpha$ depends only
upon the symmetric function ${\rm   sym}(m_\alpha)$.

The set of symmetric functions ${\rm   sym}(m_\alpha)$
obtained
as $\alpha$ varies over all sequences satisfying conditions
(i) and (ii)
equals ${\cal    M}_{\ssym}(d,n)$.
The matrix of linear equations (\ref{fredqq}) with
rows indexed by ${\cal    M}_{\ssym}(d,n)$
and columns indexed by the variable set ${\cal    P}(d,n)$
differs from the matrix $\matb^{\Bbb S}$ defined
in Section 2.5 for $(r,s)=(d,n)$
only by scalar column operations.
By Lemma \ref{lmbbs}, $\matb^{{\Bbb S}}$ has rank
equal to $|{\cal    P}(d,n)|$.
Hence the linear equations (\ref{fredqq}) uniquely
determine the integrals of partition length $n$ 
in terms of $\lan 1, \ldots, 1 \ran_g$.
The proof of the induction step is complete.
\enddemo

Since we have already found a nontrivial solution
(\ref{wqw}) of the degree $d$ localization relations
(\ref{fred}), we may conclude all solutions are
proportional to (\ref{wqw}). By Proposition \ref{cbas},
the $\lambda_g$ conjecture is proven.

\section{Partition matrices A--E}
\advance\eqcount by 26

3.1. {\it Proof of Lemma} \ref{lmaa}.
Let $r\geq s>0$.
Let ${\bf  A}$ be the matrix with rows and
columns indexed by $\pord(r,s)$ and matrix elements:
$$\mata(X,Y)= m_{Y^-}(x_1, \ldots, x_s),$$
as defined in Section 2.2.
We will prove that the matrix $\mata$ is invertible.

The set $\pord(r,s)$ may be viewed as a subset of points of
 $\com^s$ (see \S 2.2).
Matrix $\mata$ is invertible if and only if these
points impose independent conditions on the
space ${\rm   Sym}^{r-s}(\com^s)^*$
of homogeneous polynomials of degree
$r-s$ in the variables $z_1, \ldots, z_s$.

Let $v=(v_1, \ldots,v_s)$ be $s$ independent vectors
in $\com^s$.
Let $\pord(r,v)$ denote the
set of
points
$$\left\{ \ \sum_{i=1}^s x_i v_i \ | \ X=(x_1, \ldots,x_s)\in 
\pord(r,s)\ \right\}.$$
If $v$ is the standard coordinate basis,
the set $\pord(r,v)$
is the usual embedding of $\pord(r,s)$ in $\com^s$.
We will prove $\pord(r,v)$ imposes
independent conditions on ${\rm   Sym}^{r-s}(\com^s)^*$
for any basis $v$.

If $s=1$, then the cardinality of $\pord(r,s)$ is
1. The point $rv_1\neq 0$ clearly imposes
a nontrivial condition on ${\rm   Sym}^{r-1}(\com)^*$.

Let $s>1$.
By induction, we may assume   $\pord(r',v=(v_1,\ldots,v_{s'}))$
imposes independent conditions on ${\rm   Sym}^{r'-s'}(\com^{s'})^*$
for pairs $(r',s')$ satisfying $s'<s$.

If $r=s$, then the cardinality of $\pord(s,s)$
is again 1. The point $\sum_{i=1}^s v_i$
imposes a nontrivial condition on ${\rm   Sym}^0(\com^s)^*$.

Let $r>s$.
By induction, we  
may  assume $\pord(r',v=(v_1,\ldots,v_{s}))$
imposes independent conditions on ${\rm   Sym}^{r'-s}(\com^{s})^*$
for pairs $(r',s)$ satisfying $r'<r$.

We must now prove 
the points $\pord(r,v)$ impose independent
conditions on ${\rm   Sym}^{r-s}(\com^s)^*$
for any set of
independent vectors $v=(v_1, \ldots,v_s)$.
Let $f(z)\in {\rm   Sym}^{r-s}(\com^s)^*$
satisfy:
$f(p)=0$ for all $p\in \pord(r,v)$.
It suffices to prove $f(z)=0$.

Fix $1\leq j\leq s$.
Consider first the subset
\begin{equation}
\label{sss}
\left\{ \ \sum_{i=1}^s x_i v_i \ | \ X=(x_1, \ldots,x_s)\in 
\pord(r,s), \ x_j=1\ \right\} \subset \pord(r,v).
\end{equation}
The points (\ref{sss}) span a linear
subspace $L_j$ of dimension $s-1$ in $\com^s$.
In fact, the set (\ref{sss}) equals the set:
\begin{equation}
\label{sssd}
\left\{ \ \sum_{i\neq j} x_i \tilde{v}_i \ | \ \hat{X}
=(x_1, \ldots,\hat{x}_j,
\ldots,x_{s})\in 
\pord(r-1,s-1)\ \right\}
\end{equation}
where the vectors
$$\tilde{v}_i= v_i + \frac{1}{r-1} \cdot v_j, \ \ i\neq j$$
span a basis of $L_j$.
The restriction $f|_{L_j}$ lies in ${\rm   Sym}^{r-s}(L_j)^*$
and vanishes at the points (\ref{sssd}).
By our induction assumption on $s$, the restriction of
$f$ to $L_j$ vanishes identically.

The distinct linear equations defining
$L_1, \ldots, L_s$ must therefore divide $f$:
$$f= f' \cdot \prod_{i=1}^s (L_i),$$
where $f'\in {\rm   Sym}^{r-2s}(\com^s)^*$.
If $r<2s$, we conclude $f=0$.

We may assume $r\geq 2s$.
The product $\prod_{i=1}^s (L_i)$ does not
vanish at any point in the subset
\begin{equation}
\label{sssq}
\left\{ \ \sum_{i=1}^s x_i v_i \ | \ X=(x_1, \ldots,x_s)\in 
\pord(r,s), \   x_i\geq 2 \ {\rm   for all} \ i \right\} \subset \pord(r,v).\hskip.3in
\end{equation}
Hence, $f'$ must vanish at every point of (\ref{sssq}).

Define new vectors $\tilde{v}=(\tilde{v}_1, \ldots,
\tilde{v_s})$ of $\com^s$ by
$$\tilde{v}_i= \sum_{j=1}^s v_j (\delta_{ij}+\frac{1}{r-s}).$$
A straightforward determinant calculation shows
$\tilde{v}$ spans a basis of $\com^s$.
The set (\ref{sssq}) equals the set:
\begin{equation}
\label{sssqq}
\left\{ \ \sum_{i=1}^s x_i \tilde{v}_i \ | \ X=(x_1, \
\ldots,x_{s})\in 
\pord(r-s,s)\ \right\}.
\end{equation}
By the induction assumption on $r$, the function
$f'$ must vanish identically. 
We have thus proven $f=0$.
\enddemo

D. Zagier has provided us with another proof of Lemma \ref{lmaa}
by an explicit computation of the determinant:
$$|\det({\bf  A})|=
r^{\binom{r-1}{s}} \prod_{X\in\pord(r,s)} x_1^{r-s+1-x_1}.$$
We omit the derivation.

\demo{{\rm 3.2.} Proof of Lemma {\rm \ref{bbbb}}}
Let $r>s>0$.
Let $\matb$ be the matrix with rows indexed by
${\cal    M}(r,s)$, columns indexed by $\pord(r,s)$,
and matrix elements:
$$\matb(m,X)= m\left(\frac{1}{x_1}, \ldots, \frac{1}{x_s}\right),$$
as defined in Section 2.5.
We will prove matrix $\matb$ has rank
equal to $|\pord(r,s)|$.

Consider first the case $s=1$.
The set $\pord(r,1)$ consists of a single element $(r)$.
As $r\geq 2$, the constant monomial $1$
lies in $ {\cal    M}(r,1)$.
Hence $\matb$ certainly has rank equal to 1 in this case.

We now proceed by induction on $s$.
Let $s\geq 2$. 
Assume Lemma \ref{bbbb} is true for all pairs $(r',s')$
satisfying $s'<s$.

There is a natural inclusion of sets
$$\pord(r-1,s-1) \hookrightarrow \pord(r,s)$$
defined by:
$$(x_1, \ldots,x_{s-1}) \rarr (x_1,\ldots, x_{s-1}, 1).$$
Let $\pord(r-1,s-1,1)$ denote the image of this
inclusion.

There is a natural
inclusion of sets
$${\cal    M}(r-1,s-1) \hookrightarrow {\cal    M}(r,s)$$
obtained by multiplication by $z_s$:
$$ m(z_1,\ldots,z_{s-1}) \rarr
   m(z_1,\ldots,z_{s-1})\cdot z_s.$$
Let ${\cal    M}(r-1,s-1) \cdot z_s$
denote the image of this inclusion.

The submatrix of $\matb$ corresponding
to the rows ${\cal    M}(r-1,s-1)\cdot z_s$ and
columns $\pord(r-1,s-1,1)$ 
equals the matrix $\matb_{r-1,s-1}$ for the pair $(r-1,s-1)$.
By the induction assumption, we conclude
the submatrix of
columns of $\matb$ corresponding to $\pord(r-1,s-1,1)$
has full rank equal to $|\pord(r-1,s-1,1)|$.

There is a natural inclusion of sets
\begin{equation}
\label{dldl}
\pord(r-1,s) \hookrightarrow \pord(r,s)
\end{equation}
defined by:
$$(x_1, \ldots,x_{s}) \rarr (x_1,\ldots, x_{s-1}, 1+x_s).$$
Let $\pord(r-1,s^+)$ denote the image of this
inclusion.
$\pord(r,s)$ is the disjoint union of
$\pord(r-1,s-1,1)$ and $\pord(r-1,s^+)$.
We now study the columns of $\matb$ corresponding to
$\pord(r-1,s^+)$.

Let $T(z_1,\ldots,z_s)$ denote the 
polynomial function:
$$T(z)=  \left( \sum_{i=1}^{s-1} \frac{1}{z_i} \ - \ \frac{r-1}{z_s}\right) \cdot
\prod_{i=1}^s z_i.$$ \pagebreak

\specialnumber{5}\proclaim{Proposition} The function
$T(z)$ has the following properties\/{\rm :}
\begin{itemize}
\item[{\rm (i)}] $T(z)$ is homogeneous of degree $s-1$.
\item[{\rm (ii)}] Let $X\in \pord(r,s)$. Then{\rm ,}
$$T\left(\frac{1}{x_1}, \ldots, \frac{1}{x_s}\right)=0 
\ \ \leftrightarrow \ \ X\in 
\pord(r-1,s-1,1).$$
\item[{\rm (iii)}] Let $f(z)$ be any {\rm (}\/possibly nonhomogeneous\/{\rm )}
polynomial function  
of degree at most $r-s-1$. Then{\rm ,}
$$f\cdot T(z)$$
is a linear combination of monomials in ${\cal    M}(r,s)$.
\end{itemize}

\endproclaim

\demo{Proof}
Property (i) is clear by definition.
For $X\in \pord(r,s)$,
$$\sum_{i=1}^{s-1} \frac{1}{1/x_i}= r-x_s.$$
Hence $T(1/x_1, \ldots, 1/x_s)=0$ if and only if
\begin{equation}
\label{ptp}
r-x_s = (r-1) x_s.
\end{equation}
Equation (\ref{ptp}) holds if and only if $x_s=1$. Property (ii)
is thus proven.
Certainly the polynomial $f\cdot T(z)$ is of degree
at most $r-2$.
Note each monomial in $T(z)$ omits exactly 1 coordinate
factor. 
Hence each monomial of $f\cdot T(z)$ may omit at most
1 coordinate factor.
Property (iii) then holds by the definition of ${\cal    M}(r,s)$.
\enddemo

Let $\com_{r-s-1}[z]$ be the vector
space of all polynomials of 
degree at most
$r-s-1$ in  the variables $z_1, \ldots,z_s$.
Let $\com_{r-s-1}[z]\cdot T$ be the vector space of functions
$$\{ \ f\cdot T \ | \ f\in \com_{r-s-1}[z]\ \}.$$
By property (iii) of $T$, after applying 
{\it     row} operations to $\matb$, we may take the 
first ${\rm   dim}(\com_{r-s-1}[z])$ rows 
to correspond to a basis of the function space
$\com_{r-s-1}[z]\cdot T$.
Let $\matb'$ denote the matrix $\matb$ after these
row operations. 
The ranks of the column spaces of a matrix 
do {\it     not} change after row operations. Hence,
the rank of $\matb'$ equals the rank of $\matb$.
Moreover, the rank of the column space
$\pord(r-1,s-1,1)$ of $\matb'$ remains $|\pord(r-1,s-1,1)|$.

By property (ii), the block of
$\matb'$ determined by 
the row space $\com_{r-s-1}[z]\cdot T$  
and columns set $\pord(r-1,s-1,1)$ vanishes:
\begin{equation}
\label{vbk}
\matb'[\ \com_{r-s-1}[z]\cdot T,\  \pord(r-1,s-1,1)\ ] =0.
\end{equation}
Let $M$ be the
block $\matb'[\ \com_{r-s-1}[z]\cdot T,\  \pord(r-1,s^+)\ ].$
The matrix $M$ has elements:
$$M(f\cdot T, X)= f\cdot T\ \left(\frac{1}{x_1}, 
\ldots, \frac{1}{x_s}\right).$$
Since the column space $\pord(s-1,r-1,1)$ of $\matb'$ has
rank $|\pord(r-1,s-1,1)|$ and the vanishing (\ref{vbk}) holds,
$$\rk(\matb') \geq |\pord(r-1,s-1,1)| + \rk(M).$$
To prove the lemma, we will show that
the rank of $M$ equals $|\pord(r-1,s^+)|$.

Let $\matc$ be a matrix with rows indexed
by a basis of $\com_{r-s-1}[z]$,
columns indexed by $\pord(r-1,s^+)$, and matrix elements:
$$\matc(f,X)= f\left(\frac{1}{x_1}, \ldots, \frac{1}{x_s}\right).$$
As $T(1/x_1, \ldots, 1/x_s)\neq 0$ for $X\in \pord(r-1,s^+)$,
the matrix $\matc$ differs from $M$ only by scalar
column operations.
Hence,
$$\rk(M)= \rk(\matc).$$

Matrix $\matc$ is studied in Section 3.3 below. 
$\matc$ is proven to have maximal rank
$|\pord(r-1,s^+)|$ in Lemma \ref{xxxx}
by extending $\matc$ to a nonsingular
square matrix~$\matd$.

The proof of Lemma \ref{bbbb} is complete (modulo the analysis of the
matrices $\matc$ and $\matd$ in \S 3.3). 
\enddemo
 
\demo{{\rm 3.3.} Matrices $\matc$ and $\matd$}
Let $r>s>0$.
Let $\pord(\leq r-1,s)$ denote the union:
$$\pord(\leq r-1,s)= \bigcup_{t=s}^{r-1} \pord(t,s).$$
The set $\pord(\leq r-1,s)$
may be placed in bijective correspondence with a
basis of $\com_{r-s-1}[z]$ by:
\begin{equation}
\label{klkl}
X \in \pord(\leq r-1,s) \ \ \leftrightarrow \ \ m_{X^-}(z)=
z_1^{-1+x_1}\cdots z_s^{-1+x_s}.
\end{equation}
Let $\matd$ be a matrix with rows and columns indexed by
$\pord(\leq r-1,s)$.
The matrix elements of $\matd$ are
defined by:
$$\matd(X,Y)= m_{X^-}\left(\frac{1}{y_1}, \ldots, \frac{1}{y_{s-1}},
\frac{1}{1+y_s}\right).$$
Matrix $\matd$ is invertible by the following result.\enddemo

\specialnumber{6}\proclaim{Lemma}
\label{invd}
The determinant {\rm (}\/up to sign\/{\rm )} of $\matd$ is\/{\rm :}
\vglue6pt

\centerline{${\displaystyle {|\det(\matd)|}= \prod_{X\in \pord(\leq r-1,s)} 
m_{X^-}\left(\frac{1}{x_1}, \ldots, \frac{1}{x_{s-1}},
\frac{1}{1+x_s}\right) \cdot \frac{1}{x_s}.}$}
\endproclaim

{\it Proof}.
We first introduce required terminology.
For $$A=(a_1, \ldots,a_s)\in \pord(\leq r-1,s),$$
let $|A|=\sum_{i=1}^s a_i$ be the {\it     size} of $A$.
There is a partial ordering of $\pord(\leq r-1,s)$
by size.
Choose a total ordering of $\pord(\leq r-1,s)$
which refines the size partial order
(the order within each
size class may be chosen arbitrarily). This total order
of $\pord(\leq r-1,s)$ will 
be fixed for the entire proof.

Define another
partial ordering on the set $\pord(\leq r-1,s)$
by:
$$
A\geq B 
\ \ \leftrightarrow \ \
 a_i \geq b_i \ \ {\rm   for all} \ i \in \{1,  \ldots,s\}.$$
If $A\geq B$, then
either $|A|>|B|$ or $A=B$.
Hence, $B$ cannot appear strictly after $A$ in the
total order.

Let $x_1$ and $x_2$ be integers.
Define the coefficients $e_k[x_1,x_2]$ by
$$ \prod_{j=x_1}^{x_2} (t+j) = \sum_{k=0}^{x_2-x_1+1}e_k[x_1,x_2] 
\cdot t^{x_2-x_1+1-k}.$$
Note that $e_k[x_1,x_2]$ vanishes when $k>x_2-x_1+1$.
Also, $e_k[x_1,x_2]$ vanishes when $x_1>x_2$ except 
for the case $e_0(x_1, x_1-1)=1$.

The key to the proof is the construction of a related matrix
$\matd'$ with rows and columns indexed by $\pord(\leq r-1,s)$
in the fixed total order.
The matrix elements of $\matd'$ are defined 
in the following manner:
\begin{enumerate}
\item[(i)] If $A\geq B$, then
$$\matd'(A,B)= (-1)^{|B|}
\prod_{i=1}^{s-1}
\frac{e_{b_i-1}[1,a_i-1]}{(a_i-1)!} \cdot  
\frac{e_{b_s-1}[2,a_s]}{(a_s)!}.$$
\item[(ii)] In all other cases, $\matd'(A,B)=0$.
\end{enumerate}
$\matd'$ is a lower-triangular matrix with
diagonal elements:
 $$\matd'(A,A)= (-1)^{|A|}.$$
Hence $|\det(\matd')|=1$.

We now study the product $\matd' \matd$.
Consider the matrix element 
$\matd'\matd\ (A,Y)$:
$$
\sum_{i=1}^{s} \sum_{b_i=1}^{a_i} 
(-1)^{|\sum_{i=1}^s b_i|}
\prod_{i=1}^{s-1}
\frac{e_{b_i-1}[1,a_i-1]}{(a_i-1)!} \frac{1}{y_i^{b_i-1}} \cdot  
\frac{e_{b_s-1}[2,a_s]}{(a_s)!} \frac{1}{(y_s+1)^{b_s-1}}.$$
The above expression may be written in a factorized form:
$$
(-1)^s \prod_{i=1}^{s-1} \sum_{b_i=1}^{a_i} 
\frac{(-1)^{b_i-1}
e_{b_i-1}[1,a_i-1]}{(a_i-1)! \ y_i^{b_i-1}}\ \cdot \
 \sum_{b_s=1}^{a_s} 
\frac{(-1)^{b_s-1} e_{b_s-1}[2,a_s]}{(a_s)!\ (y_s+1)^{b_s-1}}\ .$$
These factors are easily evaluated. For $1\leq i \leq s-1$,
\begin{equation}
\label{firi}
 \sum_{b_i=1}^{a_i} 
\frac{(-1)^{b_i-1}
e_{b_i-1}[1,a_i-1]}{(a_i-1)! \ y_i^{b_i-1}} =
\frac{\prod_{j=1}^{a_i-1} (y_i-j)}{(a_i-1)! \ y_i^{a_i-1}}.
\end{equation}
For $i=s$,
\begin{equation}
\label{lasi}
 \sum_{b_s=1}^{a_s} 
\frac{(-1)^{b_s-1}
e_{b_s-1}[2,a_s]}{(a_s)! \ (y_s+1)^{b_s-1}} =
\frac{\prod_{j=2}^{a_s} (y_s+1-j)}{(a_s)! \ (y_s+1)^{a_s-1}}.
\end{equation}

We claim $\matd' \matd$ is upper-triangular.
Suppose $Y$ strictly precedes $A$ in the total order.
There must be a coordinate $y_i$ which satisfies
$y_i<a_i$.
If $1\leq i \leq s-1$, then the factor (\ref{firi}) vanishes.
If $i=s$, then the factor (\ref{lasi}) vanishes.
In either case, $\matd'\matd\ (A,Y)=0$.

The diagonal elements of $\matd' \matd$
are easily calculated by equations (\ref{firi}--\ref{lasi}):
$$\matd'\matd\ (A,A)= (-1)^s \prod_{i=1}^{s-1}\frac{1}{a_i^{a_i-1}}
\cdot \frac{1}{(a_s+1)^{a_s-1} a_s}.$$
As $\matd' \matd$ is upper-triangular, the
determinant is the product of the diagonal entries.
Since $|\det(\matd')|=1$, this determinant equals (up to sign)
$\det(\matd)$.
\hfill\qed\vglue12pt

Consider the column set of $\matd$
corresponding to the subset
$$\pord(r-1,s) \subset \pord(\leq r-1,s).$$
The submatrix of $\matd$ obtained by restriction to
this column set equals
$\matc$ via the correspondences
(\ref{klkl}) and (\ref{dldl}).
As a corollary of Lemma \ref{invd}, we may conclude the
required rank result for $\matc$.

\specialnumber{7}\proclaim{Lemma}
\label{xxxx}
Let $r>s>0$.
$\matc$ has rank equal to $|\pord(r-1,s^+)|$.
\endproclaim

3.4. {\it Matrix ${\bf  E}$}.
Let ${\bf  E}$ be a matrix with row and columns
indexed by the set $\pord(\leq r-1,s)$.
The matrix elements of ${\bf  E}$ are
defined by:
$${\bf  E}(X,Y)= m_{X^-}\left(\frac{1}{y_1}, \ldots,
\frac{1}{y_s}\right).$$
While
${\bf  E}$ is slightly more natural than $\matd$,
we do not encounter ${\bf  E}$ in our
proof of the
$\lambda_g$ conjecture.
We note, however, the proof of Lemma \ref{invd}
may be modified to prove:

\specialnumber{8}\proclaim{Lemma}
The determinant {\rm (}\/up to sign\/{\rm )} of ${\bf  E}$ is\/{\rm :}
$${|\det({\bf  E})|}= \prod_{X\in \pord(\leq r-1,s)} 
m_{X^-}\left(\frac{1}{x_1}, \ldots, 
\frac{1}{x_s}\right) .$$
\endproclaim
 
\AuthorRefNames [COGP]

\end{document}